\documentclass[12pt]{article}
\usepackage{mathrsfs}
\usepackage{latexsym,amsmath,amssymb,amsfonts,epsfig,graphicx,cite}
\newtheorem{theorem}{Theorem}[section]

\newtheorem{proposition}[theorem]{Proposition}

\newtheorem{prop}[theorem]{Property}

 \numberwithin{equation}{section}
 \numberwithin{figure}{section}

\def\qed{\hfill \rule{4pt}{7pt}}

\def\pf{\noindent {\it Proof.} }

\begin{document}

\begin{center}
{\large\bf On the Combinatorics of the Boros-Moll Polynomials}
\end{center}

\vskip 2mm \centerline{William Y.C. Chen$^1$, Sabrina X.M. Pang$^2$,
and Ellen X.Y. Qu$^3$}

\begin{center}
$^{1, 2, 3}$Center for Combinatorics, LPMC-TJKLC\\
Nankai University, Tianjin 300071, P.R. China

\vskip 2mm
  $^1$chen@nankai.edu.cn,  $^2$pang@cfc.nankai.edu.cn,
 $^3$xiaoying@cfc.nankai.edu.cn
\end{center}

\begin{center}
{\bf Abstract}
\end{center}

The Boros-Moll polynomials arise in the evaluation of a quartic
integral. The original double summation formula does not imply the
fact that the coefficients of these polynomials are positive.  Boros
and Moll proved the positivity by using Ramanujan's Master Theorem
to reduce the double sum to a single sum. Based on the structure of
reluctant functions introduced by Mullin and Rota along with an
extension of Foata's bijection between Meixner endofunctions and
bi-colored permutations, we find a combinatorial proof of the
positivity. In fact, from our combinatorial argument one sees that
it is essentially the binomial theorem that makes it possible to
reduce the double sum to a single sum.

\vskip 3mm


\vskip 3mm



\section{Introduction}

The objective of this paper is to give a combinatorial proof of the
positivity of the coefficients of the Boros-Moll polynomials.
 Boros and Moll
\cite{Bor-Mol99a,Bor-Mol99b,Bor-Mol99c,Bor-Mol01,Bor-Mol04,Mol02}
explored the following integral which is closely related to a
special class of Jacobi polynomials. They have shown that for any
$a>-1$ and any nonnegative integer $m$,
\begin{equation}
\int_{0}^{\infty}\frac{1}{(x^4+2ax^2+1)^{m+1}}dx=\frac{\pi}{2^{m+3/2}(a+1)^{m+1/2}}P_m(a),
\end{equation}
where
\begin{equation}\label{eq1.2}
P_m(a)=\sum_{j, k}{2m+1\choose 2j}{m-j\choose k}{2k+2j\choose
k+j}\frac{(a+1)^j(a-1)^k}{2^{3(k+j)}}.
\end{equation}
The polynomials $P_m(a)$ are called the Boros-Moll polynomials
\cite{Chen-Xia08}.  Write
$$P_m(a)=\sum\limits_{i=0}^{m}d_i(m)a^i.$$
Boros and Moll found a remarkable proof of the fact that the
coefficients $d_i(m)$ are positive by employing  Ramanujan's Master
Theorem, see \cite{Bor-Mol01} or \cite[Theorem 7.9.1]{Bor-Mol04}. In
fact, they have  shown that
\begin{equation}\label{eq1.3}
P_m(a)=2^{-2m}\sum\limits_{k}2^k{2m-2k\choose m-k}{m+k\choose
k}(a+1)^k.
\end{equation}
It follows from \eqref{eq1.3} that
$$d_i(m)=2^{-2m}\sum\limits_{k=i}^m 2^k{2m-2k\choose m-k}{m+k\choose
k}{k\choose i}.$$ There are several  proofs of this formula, see the
survey of  Amdeberhan and  Moll \cite{Amd-Mol08}. By the above
formula (\ref{eq1.3}), one can express   $P_m(a)$ in terms of a
hypergeometric series
$$P_m(a)=2^{-2m}{2m\choose m}{_2F_1\left(-m, m+1; \frac{1}{2}-m; \frac{a+1}{2}\right)}.$$
Recall that ${_2F_1}$ denotes the hypergeometric series
$${_2F_1(a, b; c; x)}=\sum\limits_k\frac{(a)_k(b)_k}{(c)_k}\cdot
\frac{x^k}{k!},$$ where
 $(a)_k$ stands for the rising factorial defined by
  $(a)_k=a(a+1)\cdots(a+k-1)$ for $k>0$ and
$(a)_k=1$ for $k=0$. Consequently, $P_m(a)$ can be viewed as the
Jacobi polynomial $P_m^{(\alpha, \beta)}(a)$ with
\[ \alpha=m+\frac{1}{2}, \quad
\beta=-m-\frac{1}{2}.\] Recall that $P_m^{(\alpha, \beta)}(a)$ is
defined by
$$P_m^{(\alpha, \beta)}(a)=\sum\limits_{k=0}^m(-1)^{m-k}{m+\beta\choose m-k}
{m+k+\alpha+\beta\choose k}\left(\frac{1+a}{2}\right)^k.$$

Much progress has been made since Boros and Moll proved the
positivity of the coefficients of $P_m(a)$. Boros and Moll
\cite{Bor-Mol99c} have shown that the sequence $\{d_i(m)\}_{0\leq i
\leq m}$ is unimodal, that is, there exists an index $i$ such that
$d_{0}(m)\leq \cdots \leq d_{i}(m)$ and $d_{i}(m)\geq \cdots \geq
d_{m}(m)$. Moll conjectured that the coefficients $d_i(m)$ form a
log-concave sequence, that is, $ d_{i}(m)^2\geq
d_{i-1}(m)d_{i+1}(m)$ for $1\leq i \leq m-1$. This conjecture has
been confirmed by Kauers and Paule \cite{Kau-Pau07}. Recently, Chen
and Xia \cite{Chen-Xia08} have proved a stronger property of
$d_i(m)$, called the ratio monotone property, which implies the
log-concavity and the spiral property.  The combinatorial aspects of
the $2$-adic valuation of the number $i!m!2^{m+i}d_i(m)$ have been
studied by Amdeberhan, Manna and Moll \cite{Amd-Man-Mol08}, and Sun
and Moll \cite{Sun-Mol08}.

From the combinatorial point of view,  it is always interesting to
find combinatorial reasons for the coefficients to be positive when
the direct expansion contains negative terms. It is also desirable
to find combinatorial proofs of unimodal and log-concave properties.
Furthermore, it would be interesting to find combinatorial
interpretations of the recurrence relations of $d_i(m)$ given by
Kauers and Paule \cite{Kau-Pau07} and Moll \cite{Mol07}.

In this paper, we will take the first step in this direction. We
will give a combinatorial interpretation of the positivity of the
coefficients of the Boros-Moll polynomials based on the structure of
reluctant functions introduced by Mullin and Rota \cite{Mul-Rot95}
along with an extension of Foata's bijection between Meixner
endofunctions and bi-colored permutations. It should noted that the
structure of reluctant functions and Meixner endofuntions have also
been used in the combinatorial study of the Pfaff identity by Chen
and Pang \cite{Chen-Pang08}.

More specifically, we will give a combinatorial proof of the
following identity which implies the equivalence of the two
expressions \eqref{eq1.2} and \eqref{eq1.3} for $P_m(a)$:
\begin{eqnarray}
\lefteqn{\sum_{j,k} {2m+1\choose 2j}{m-j\choose k}{2k+2j\choose
k+j}\frac{(a+1)^j(a-1)^k}{2^{3(k+j)}}} \qquad \nonumber\\[3pt]
&=& 2^{-2m}\sum_{k}2^k{2m-2k\choose m-k}{m+k\choose k}(a+1)^k.
\label{21}
\end{eqnarray}

\section{The Combinatorial Proof.}

 In order to give a combinatorial interpretation of the relation
 (\ref{21}) that implies the positivity of the coefficients of the Boros-Moll identities,
 we need to use a variant of the identity by multiplying both sides by $m!$.
 The following reformulation of the identity after the multiplication by $m!$
 is straightforward and can be made
 purely combinatorial in principle.
  Let us
denote the left hand side and the right hand side of \eqref{21} by
$L$ and $R$, respectively. Throughout this paper, we will use the
notation $(x)_n$ to for rising factorials, that is,
$(x)_n=x(x+1)\cdots(x+n-1)$ for $n>0$ and $(x)_n=1$ for $n=0$. On
one hand, we have \allowdisplaybreaks
\begin{align*}
m!\cdot L&=m!\sum_{i+j+k=m}\frac{(2m+1)!}{(2j)!(2m+1-2j)!} \cdot
\frac{(m-j)!}{k!i!}\cdot \frac{(2m-2i)!}{(m-i)!(m-i)!}\cdot
\frac{(a+1)^j(a-1)^k}{2^{3m-3i}}\\
&=m!\sum_{i+j+k=m}\frac{2^{2m+1}m!(m+\frac{1}{2})!}
{2^{2j}j!(j-\frac{1}{2})!2^{2m+1-2j}(m+\frac{1}{2}-j)!(m-j)!} \cdot
\frac{(m-j)!}{k!i!}\\
&\hskip 2cm \cdot
\frac{2^{2m-2i}(m-i)!(m-i-\frac{1}{2})!}{(m-i)!(m-i)!}\cdot
\frac{(a+1)^j(a-1)^k}{2^{3m-3i}}\\
&=\sum_{i+j+k=m}{m\choose i,j,k}\frac{m!(m+\frac{1}{2})!}
{(j-\frac{1}{2})!(m+\frac{1}{2}-j)!} \cdot
\frac{(m-i-\frac{1}{2})!}{(m-i)!}\cdot
\left(\frac{a+1}{2}\right)^j\left(\frac{a-1}{2}\right)^k\\
&=\sum_{i+j+k=m}{m\choose
i,j,k}(m-i+1)_i\left(m-j+\frac{3}{2}\right)_j\left(j+\frac{1}{2}\right)_k\left(\frac{a+1}{2}\right)^j
\left(\frac{a-1}{2}\right)^k\\
&=\sum_{i+j+k=m}{m\choose
i,j,k}(-m)_i(-1)^i\left(-m-\frac{1}{2}\right)_j(-1)^j
\left(j+\frac{1}{2}\right)_k\left(\frac{a+1}{2}\right)^j
\left(\frac{-a+1}{2}\right)^k(-1)^k\\
&=(-1)^m\sum_{i+j+k=m}{m\choose
i,j,k}(-m)_i\left(-m-\frac{1}{2}\right)_j
  \left(j+\frac{1}{2}\right)_k\left(\frac{a+1}{2}\right)^j
\left(\frac{-a+1}{2}\right)^k.
\end{align*}
On the other hand, we have
\begin{align*}
m!\cdot R&=m!\cdot\sum_{i+j=m}2^{-2m+i}{2j\choose
j}{m+i\choose i}(a+1)^i\\
&=m!\cdot\sum_{i+j=m}2^{-2m+i}\frac{(2j)!}{j!j!}
\frac{(m+i)!}{i!m!}(a+1)^i\\
&=\sum_{i+j=m}{m\choose i,j}2^{-2m+i}\frac{(2j)!(m+i)!}
{j!m!}(a+1)^i\\
&=\sum_{i+j=m}{m\choose i,j}2^{-i}\left(j-\frac{1}{2}\right)!
\frac{(m+i)!}
{m!}(a+1)^i\\
&=\sum_{i+j=m}{m\choose
i,j}2^{-i}(-1)^j\left(\frac{1}{2}-j\right)_j(m+1)_i
(a+1)^i\\
&=(-1)^m\sum_{i+j=m}{m\choose
i,j}\left(\frac{1}{2}-j\right)_j(m+1)_i
\left(\frac{-a-1}{2}\right)^i.
\end{align*}
So the identity \eqref{21} can be converted into the following
equivalent form
\begin{eqnarray}\label{eq2.1}
\lefteqn{\sum_{i+j+k=m}{m\choose
i,j,k}(-m)_i\left(-m-\frac{1}{2}\right)_j
  \left(j+\frac{1}{2}\right)_k\left(\frac{a+1}{2}\right)^j
\left(\frac{-a+1}{2}\right)^k}  \qquad \qquad \qquad\nonumber\\
& =&\sum_{i+j=m}{m\choose i,j}\left(\frac{1}{2}-j\right)_j(m+1)_i
\left(\frac{-a-1}{2}\right)^i. \qquad \qquad \qquad
\end{eqnarray}

Our combinatorial approach to the above identity consists of three
steps. The first step is to give combinatorial interpretations of
the sums on both sides of \eqref{eq2.1}. We will show that the left
hand side is the sum of weights of Meixner bi-endofunctions, and the
right hand side is the sum of weights of Meixner endofunctions with
a different weight assignments. The second step is to transform the
sum of weights of Meixner bi-endofunctions to the sum of weights of
$3$-colored permutations. This is achieved by a weight preserving
bijection  between Meixner bi-endofunctions and 3-colored
permutations, which is a natural extension of Foata's bijection.
Meanwhile, the sum of weights for the right hand side can be
transformed to the sum of weights of bi-colored permutations by the
original bijection of Foata.  The third step is to compare the
weights of $3$-colored permutations and bi-colored permutations. One
sees that the equality follows from the weight distribution on a
cycle. Roughly speaking, if there are two ways to give a weight
$w_1$ or $w_2$ to a cycle, then it is equivalent to assigning only
one weight $w_1+w_2$ to the cycle. This step yields a combinatorial
interpretation of why the double sum (\ref{eq1.2}) reduces to a
single sum (\ref{eq1.3}).

Note that a basic ingredient of the combinatorial settings for the
above hypergeometric identity is the interpretation of the rising
factorial $(x)_n$, or, in general, of $(x+k)_n$. It is well known
that $(x)_n$ can be expanded in terms of the signless Stirling
numbers of the first kind. Note that $(x)_n$ can also be interpreted
as the number of dispositions from $[n]=\{1, 2, \ldots, n\}$ to a
set $X$ with $x$ elements, see Joni, Rota and Sagan
\cite{Jon-Rot-Sag81} for more details.

In general, the rising factorial $(a+j)_i$ can be explained as the
sum of  the weights of reluctant functions from $A$ to $B$, where
$A$ and $B$ are disjoint, and $|A|=i$ and $|B|=j$. Recall that the
notion of reluctant functions was introduced by Mullin and Rota
\cite{Mul-Rot95} in their theory of sequences of polynomials of
binomial type. A reluctant function $f$ from $A$ to $B$, where $A$
and $B$ are two disjoint finite sets, is defined as an injective map
from $A$ to $A\cup B$. The functional digraph of $f$ is a digraph on
$A\cup B$ with arcs $(k, f(k))$ for $k\in A$. The weight of $f$ is
defined as $a^k$, where $k$ is the number of cycles in the
functional digraph of $f$.

Observe that the functional digraph of any reluctant function $f$
has a unique decomposition into disjoint cycles on elements in $A$
and directed paths ending with an element in $B$. The ending points
in $B$ are called {\it terminals}. Now, let us review the canonical
cycle representation of a reluctant function,  introduced by Chen
and Pang in \cite{Chen-Pang08} as a natural extension of the
canonical cycle representation of a permutation, see Stanley
\cite[Page 17]{Sta99}. Assume that $f$ is a reluctant function from
$A$ to $B$. The functional digraph of $f$ can be decomposed into $k$
cycles $C_1, C_2, \ldots, C_k$ and $s$ directed paths $P_1, P_2,
\ldots, P_s$. We first write down the cycles in canonical cycle
representation, that is, write a cycle $C=(i_1i_2\cdots i_r)$  in
such a way that  $i_1$ is the minimum element of $C$,  then arrange
the cycles $C_1, C_2, \ldots, C_k$  in accordance with the
decreasing order of their minimum elements. Moreover, each path
$P_i$ is written as $(j_1j_2\cdots j_l)$ such that $j_1\in B$ and
$f(j_t)=j_{t-1}$ for $2\leq t\leq l$, and $P_1, P_2, \ldots, P_s$
are arranged according to the increasing order of  their first
elements.

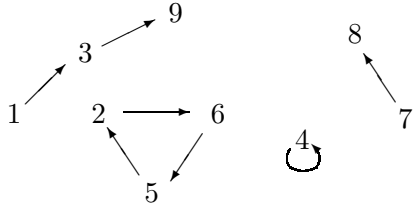
\begin{figure}
\begin{center}
\begin{picture}(200,50)(0,0)
\put(-2,17){\small$1$}\put(5,24){\vector(1,1){15}}
\put(25,40){\small$3$}\put(34,46){\vector(2,1){20}}
\put(59,55){\small$9$} \put(30,17){\small$2$}
\put(42,21){\vector(1,0){25}} \put(75,17){\small$6$}
\put(50,-13){\small$5$}
\put(48,-3){\vector(-2,3){12}}\put(72,13){\vector(-2,-3){12}}
\put(145,25){\vector(-2,3){12}}\put(146,15){\small$7$}
\put(127,47){\small$8$}\put(107,7){\small$4$}
\qbezier(105,7)(102.5,3.5)(105,0)\qbezier(105,0)(110,-2.5)(115,0)
\qbezier(115,0)(117.5,3.5)(115,7)\put(115,7){\vector(-1,1){2}}
\end{picture}
\end{center}\vspace{5pt}\caption{The digraph of a reluctant function.}\label{fi1}
\end{figure}

For example,  the reluctant function in Figure \ref{fi1} with
$A=\{1,2,3,4,5,6,7\}$ and $B=\{8,9\}$ has the following canonical
cycle representation
\[  (4)(265)({\bf 8}7)({\bf 9}31). \]
It can be seen that the canonical cycle representation is in fact
uniquely determined by the sequence $4 2 65 {\bf 8} 7 {\bf 9} 3 1$.
Clearly, the reluctant function $f$ can be
 recovered from the canonical cycle representation.
 To transform a sequence $a_1a_2\cdots a_m$ to the canonical cycle representation, we
 need to consider the  left-to-right minimum
elements in the sequence. Recall that an element $a_i$ in
$a_1a_2\cdots a_m$ is called a left-to-right minimum element if
$a_i<a_j$ for any $j<i$. For example, $4$ and $2$ are left-to-right
minimum elements in the above sequence.  On one hand, we can insert
a left parenthesis in the sequence before each element in $B$, which
is in boldface. On the other hand, we can insert a left parenthesis
in the sequence preceding every left-to-right minimum element $a_i$
in $A$ as long as $a_1, a_2, \ldots, a_{i-1}$ all belong to $A$.
After the left parentheses are placed in the sequence, the right
parentheses can be added accordingly.

The following proposition  is well-known,  see, for example,
 \cite{Chen-Pang08, Foa83, Foa-Str84, Lab-Yeh89, Ler-Str85}. It plays a crucial role in
the combinatorial interpretation of the identity \eqref{eq2.1}.

\begin{proposition}\label{pro2.1}
Let $A$ and $B$ be two disjoint subsets of $[m]$, and let $|A|=i$
and $|B|=j$. Then the sum of weights of reluctant functions from $A$
to $B$ equals $(a+j)_i$.
\end{proposition}

To show that the left hand side of the identity \eqref{eq2.1} equals
the sum of weights of $3$-colored permutations, we need an extension
of Foata's bijection between Meixner endofunctions and bi-colored
permutations \cite{DF,LY1}. To be more specific, we will extend
Foata's bijection to  Meixner bi-endofunctions and $3$-colored
permutations. Recall that a {\it Meixner endofunction} on a finite
set $S$ is represented by
 $(A,B;\pi_A;\pi_B)$,
where $(A,B)$ is a composition of $S$ and $\pi_A$ is an injective
map from $A$ to $S$ and $\pi_B$ is a permutation on $B$. A {\it
bi-colored permutation} on a finite set $S$ is represented by
$(A,B;\sigma)$, where $(A,B)$ is a composition of $S$, and $\sigma$
is a permutation on $S$. Note that a composition $(A, B)$ of a set
$S$ can be considered as a $2$-coloring of $S$. Foata's bijection
can be described as follows.

\begin{prop}\label{prop2.2}
There is a  bijection between the set of Meixner endofunctions on
$[m]$ and the set of bi-colored permutations on $[m]$.
\end{prop}

We now define Meixner bi-endofunctions and 3-colored permutations. A
{\it Meixner bi-endofunction} on a finite set $S$ is denoted by
$(A,B,C;\pi_A; \pi_B, \pi_C)$, where $(A,B,C)$ is a composition of
$S$, $\pi_A$ is an injective map from $A$ to $A\cup B$ and
($\pi_B$,$\pi_C$), where $\pi_B$ is a permutation on $B$ and $\pi_C$
is a permutation on $C$.
 Given a 3-coloring of $[m]$, say by the three colors
red, black and white, a  $3$-colored permutation is defined as a
permutation on $[m]$ such that no red elements appears in any cycle
containing a black or white element. For example,
$$(\underline{8},
\underline{7},\underline{9})({\bf 2},5,{\bf 4})({\bf 10},
1)({\bf3})(\underline{11},\underline{12})(6)$$ is a $3$-colored
permutation, where  the underlined elements are red, and the black
elements are in boldface.

Notice that a Meixner bi-endofunction $(A,B,C; \pi_A; \pi_B,\pi_C)$
reduces to a Meixner endofunction when $C=\emptyset$. Applying
Foata's bijection to the cycles composed of elements in $B$, we
obtain the following extension of Proposition \ref{prop2.2}.

\begin{prop}\label{prop2.3}
There is a bijection between the set of Meixner bi-endofunctions on
$[m]$ and the set of  $3$-colored permutations on $[m]$.
\end{prop}

\pf Given a Meixner bi-endofunction $(A,B,C;\pi_A;
\pi_B, \pi_C)$, we color the elements in $A,B$ and $C$
by white, black and red, respectively. Consider the
cycle representation of $\pi_A,\pi_B,\pi_C$. We may
view a Meixner bi-endofunction as a union of disjoint
cycles on $A,B,C$ along with some directed paths on $A$
attached to some element in $B$. Since $\pi_A$ is
injective, two directed paths on $A$ cannot be incident
to the same element in $B$.

The bijection will involve only the components consisting of cycles
on a subset of $B$ attached with some paths on $A$. Let $D$ be such
a cycle, and $P$ be a directed path attached to $D$. Assume that $x$
is the terminal element of $P$ that is on $D$. Let $(y,x)$ be an arc
on  on $D$. Then we can break this  arc from $y$ to $x$ and connect
$y$ to the starting point of $P$. Considering the colors of the
elements on the path $P$, we see that the above operation is
reversible. Taking all the paths attached to $D$ into account, we
obtain the desired bijection.  \qed

\begin{figure}
\begin{center}
\begin{picture}(190,50)(0,0)
\put(-2,-3){\small${\bf3}$} \put(-22,28.5){\small${\bf1}$}
\put(18,28.5){\small${\bf6}$} \put(-23,58){\small$4$}
\put(18,58){\small$5$}
\put(-7,7){\vector(-1,2){8}}\put(14,22.5){\vector(-1,-2){8}}
\put(-12,30){\vector(1,0){22}}\put(-20,52){\vector(0,-1){14}}
\put(20,52){\vector(0,-1){14}}
\put(46,73){\small$2$}\put(42,74){\vector(-2,-1){16}}
\put(70,30){$\Leftrightarrow$} \put(126,36){\vector(1,1){15}}
\put(141,9){\vector(-1,1){15}}
\put(155,58){\vector(1,0){20}}\put(175,4){\vector(-1,0){20}}
\put(190,51){\vector(1,-1){14}} \put(204,24){\vector(-1,-1){14}}
\put(182,0){\small${\bf6}$}\put(145,0){\small${\bf3}$}
\put(118,27){\small$4$}\put(145,54){\small${\bf1}$}
\put(182,54){\small$2$}\put(208,27){\small$5$}
\end{picture}
\end{center}\caption{Foata's bijection.}\label{figure6}
\end{figure}
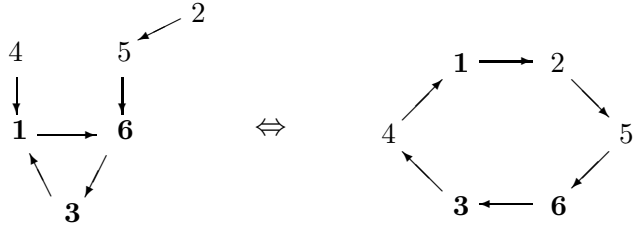

For example, as illustrated in Figure \ref{figure6}, the Meixner
endofunction
\begin{align*}
(\{2,4,5\},\{1,3,6\};(1,4)(6,5,2);(1,6,3)),
\end{align*}
corresponds to the bi-colored permutation
$({\bf3},4,{\bf1},2,5,{\bf6})$, where the black elements are in
boldface.

We are now ready to give a combinatorial proof of the identity
(\ref{eq2.1}). First, we define the weights of Meixner
bi-endofunctions and Meixner endofunctions. Let $(A,B,C; \pi_A;
\pi_B,\pi_C)$ be a Meixner bi-endofunction on $[m]$. An element in
$A,B$ or $C$ is assigned the weight
\[(-a+1)/2, \quad
(a+1)/2, \quad 1,\] respectively. Similarly, the weight of a cycle
in $\pi_A, \pi_B$ or $\pi_C$ is given by \[1/2, \quad -m-1/2, \quad
-m.\]  Then  the weight of a Meixner bi-endofunctions is the product
of the weights of the elements and the weights of the cycles.

Next, we define the weight of a Meixner endofunction $(A, B; \pi_A;
\pi_B)$ on $[m]$. The weight of an element in $A$ is given by $1$,
the weight of an element in $B$ is given by $(-a-1)/2$, the weight
of a cycle in $\pi_A$ is given by $1/2-m$, and the weight of a cycle
in $\pi_B$ is given by $1+m$. Then the weight of a Meixner
endofunction is the product of the weights of the elements and the
weights of the cycles. Given the above weight assignments for
Meixner bi-endofunctions and Meixner endofunctions,  the identity
\eqref{eq2.1} is equivalent to the following statement.

\begin{theorem} The sum of weights of Meixner
bi-endofunctions on $[m]$ equals the sum of weights of Meixner
endofunctions on $[m]$.
\end{theorem}

 \pf
By Proposition \ref{pro2.1}, it is not check to see that the sum of
weights of Meixner bi-endofunctions $(A, B,C; \pi_A; \pi_B,\pi_C)$
on $[m]$ equals the summation on the left hand side of
(\ref{eq2.1}):
\begin{eqnarray}\label{e6}
\sum_{i+j+k=m}{m\choose
i,j,k}(-m)_i\left(-m-\frac{1}{2}\right)_j
   \left(j+\frac{1}{2}\right)_k\left(\frac{a+1}{2}\right)^j
\left(\frac{-a+1}{2}\right)^k.
\end{eqnarray}

Applying the bijection  described in Proposition \ref{prop2.3}
between Meixner bi-endofunctions on $[m]$ and 3-colored permutations
on $[m]$, we find that (\ref{e6}) can be rewritten as the summation
of weights of $3$-colored permutations on $[m]$ with the following
weight assignments. A white, black, or red element is given the
weight
\[(-a+1)/2, \quad (a+1)/2, \quad 1.\] A cycle
containing only white elements is given the weight $1/2$, a cycle
containing at least one black element is given the weight $-m-1/2$,
and a cycle consisting of only red elements is given the weight
$-m$. Now, the weight of a  $3$-colored permutation is defined as
the product of the weights of the elements and the weights of the
cycles.

On the other hand, the total weight of  $3$-colored permutations on
$[m]$ can be computed based on the cycle decompositions of
permutations on $[m]$. Given a permutation $\pi$ on $[m]$ and a
cycle $D$ in $\pi$ with $r$ elements, if $D$ is a cycle consisting
of white elements, then the weight contribution is
\begin{equation}\label{m0}
\frac{1}{2}\left(\frac{-a+1}{2}\right)^r.
\end{equation}
If $D$ is used to form a cycle containing at least one black
element, the total weight contribution equals
\begin{equation} \label{m1}
\left(-m-\frac{1}{2}\right)\sum_{i=1}^r{r\choose
i}\left(\frac{a+1}{2}\right)^i \left(\frac{-a+1}{2}\right)^{r-i}.
\end{equation}
If $D$ is used to a cycle containing only red elements, the total
weight contribution equals $-m$. Combining the above three cases, we
get the total weight contribution of the cycle $D$ to the summation
of weights of  $3$-colored permutations
\[
-m+\frac{1}{2}\left(\frac{-a+1}{2}\right)^r+\left(-m-\frac{1}{2}\right)
\sum_{i=1}^r{r\choose i}\left(\frac{a+1}{2}\right)^i
\left(\frac{-a+1}{2}\right)^{r-i}  ,\] which simplifies to
\begin{equation}
-2m-\frac{1}{2}+(m+1)\left(\frac{-a+1}{2}\right)^r.\label{eq2}
\end{equation}
Note that we can easily give a combinatorial argument for the above
computation.

We continue to show that the right hand side of \eqref{eq2.1} can
also be expressed as a summation over permutations on $[m]$ with
each cycle having the above weight \eqref{eq2}. By the definition of
the weight of a Meixner endofunction, it is easily seen that the sum
of  weights over Meixner endofunctions on $[m]$ equals the the
summation on the right hand side of \eqref{eq2.1}:
\begin{equation} \label{equation7}
\sum_{i+j=m}{m\choose i,j}\left(\frac{1}{2}-j\right)_j(m+1)_i
\left(\frac{-a-1}{2}\right)^i.
\end{equation}

Applying the bijection in Proposition \ref{prop2.2} between Meixner
endofunctions on $[m]$ and bi-colored permutations on $[m]$,
\eqref{equation7} can be expressed as a summation of weights of
bi-colored permutations on $[m]$ with the following weight
assignments. The weight of a white element is given by $1$, the
weight of a black element is given by $(-a-1)/2$, the weight of a
cycle consisting of only white elements is given by $1/2-m$ and the
weight of a cycle containing at least one black element is given by
$1+m$.

Analogously, the total weight of the bi-colored permutations on
$[m]$ can be computed based on the cycle decompositions of
permutations on $[m]$. Given a permutation $\pi$ on $[m]$ and a
cycle $D$ in $\pi$ with $r$ elements, if $D$ is a cycle consisting
of white elements, the weight contribution is $1/2-m$. If $D$ is
used to form a cycle containing at least one black element, the
total weight contribution equals
\begin{equation} \label{mm1}
(1+m)\sum_{i=1}^r{r\choose i}\left(\frac{-a-1}{2}\right)^i 1^{r-i}=
(1+m)\left[\left(\frac{-a+1}{2}\right)^r-1\right].
\end{equation}
Summing up the above two cases, we get the total weight contribution
of the cycle $D$ to the summation of weights of bi-colored
permutations on $[m]$:
\begin{eqnarray}
\frac{1}{2}-m+(1+m)\left[\left(\frac{-a+1}{2}\right)^r-1\right]
=-2m-\frac{1}{2}+(m+1)\left(\frac{-a+1}{2}\right)^r.\label{eq3}
\end{eqnarray}
Comparing  \eqref{eq2} and \eqref{eq3}, we see that the weight
assignment to $3$-colored permutations is equivalent to the weight
assignment to bi-colored permutations. This completes the
combinatorial proof of the identity \eqref{eq2.1}. \qed

\vskip 3mm \noindent {\bf Acknowledgments.} This work was supported
by the 973 Project, the PCSIRT Project of the Ministry of Education,
the Ministry of Science and Technology,  the National Science
Foundation of China.


\begin{thebibliography}{99}

\bibitem{Amd-Mol08} T. Amdeberhan and V.H. Moll, A formula for a quartic integral: a
survey of old proofs and some new ones, Ramanujan J., to appear.

\bibitem{Amd-Man-Mol08} T. Amdeberhan, D. Manna, and V.H. Moll, The
$2$-adic valuation of a sequence arising from a rational integral,
J. Combin. Theory, Ser. A, to appear.


\bibitem{Bor-Mol99a} G. Boros and V.H. Moll, An ingegral hidden in
Gradshteyn and Ryzhik, J. Comput. Appl. Math. 106 (1999) 361--368.

\bibitem{Bor-Mol99b} G. Boros and V.H. Moll, A sequence of unimodal
polynomials, J. Math. Anal. Appl. 237 (1999) 272--285.

\bibitem{Bor-Mol99c} G. Boros and V.H. Moll, A critrion for
unimodality, Electron. J. Combin. 6 (1999) R3.

\bibitem{Bor-Mol01} G. Boros and V.H. Moll, The double square root,
Jacobi polynomials and Ramanujan's Master Theorem, J. Comput. Appl.
Math. 130 (2001) 337--344.

\bibitem{Bor-Mol04} G. Boros and V.H. Moll, Irresistible Integrals,
Cambridge University Press, Cambridge, 2004.


\bibitem{Chen-Pang08} W.Y.C. Chen and S.X.M. Pang, On the
combinatorics of the Pfaff identity, Discrete Math.,  to appear.

\bibitem{Chen-Xia08} W.Y.C. Chen and E.X.W. Xia, The ratio
monotonicity of the Boros-Moll polynomials, Math. Comput.,  to
appear.

\bibitem{Foa83} D. Foata, Combinatoire de identit\'{e}s sur les
polyn\^{o}mes orthogonaux, Proc. Int. Congr. Mathematicians, Warsaw,
Poland, 1983, 1541--1553.

\bibitem{DF}
D. Foata and J. Labelle, Mod\`{e}les combinatoires pour les
polyn\^{o}mes de Meixner,  European J. Combin., 4 (1983), 305-311.

\bibitem{Foa-Str84} D. Foata and V. Strehl, Combinatorics of the
Laguerre polynomials, in: Enumeration and Design, D.M. Jackson and
S.A. Vanstone, Eds., Academic Press, 1984, 123--140.

\bibitem{Jon-Rot-Sag81} S.A. Joni, G.-C. Rota and B. Sagan, From
sets to functions: Three elementary examples, Discrete Math. 37
(1981) 193--202.

\bibitem{Lab-Yeh89} J. Labelle and Y.N. Yeh, The combinatorics of
Laguerre, Charlier, and Hermite polynomials, Studies in Applied
Math. 80 (1989) 25--36.

\bibitem{LY1}
J. Labelle and Y. N. Yeh, Combinatorial proofs of some limit
formulas involving orthogonal polynomials,  Discrete Math. 79 (1990)
77--93.

\bibitem{Ler-Str85} P. Leroux and V. Strehl, Jacobi polynomials:
combinatorics of the basic identities, Discrete Math. 57 (1985)
267--187.

\bibitem{Mol02} V.H. Moll, The evaluation of integrals: A personal
story, Notices Amer. Math. Soc. 49(3) (2002) 311--317.


\bibitem{Mol07} V.H. Moll, Combinatorial sequences arising from a rational integral, Online
Journal of Analytic Combinatorics Issue 2 (2007) \#4.


\bibitem{Mul-Rot95} R. Mullin and G.-C. Rota, On the foundation of
combinatorial theory: III. Theory of binomial enumeration, Graph
Theory and Its Applications, B. Harris, ed., Academic Press, New
York and London, 1970, pp. 167--213. Reprinted in Gian-Carlo Rota on
Combinatorics, J.P.S. Kung, ed., Birkh\"{a}user,
 1995, pp. 118--147.

\bibitem{Kau-Pau07} M. Kauers and P. Paule, A computer proof of
Moll's log-concavity conjecture. Proc. Amer. Math. Soc. 135 (2007)
3837--3846.

\bibitem{Sta99} R.P. Stanley, Enumerative Combinatorics, Vol. 1,
Cambridge University Press, Cambridge, New York, 1999.

\bibitem{Sun-Mol08} X.Y. Sun and V.H. Moll, A binary tree
representation for the $2$-adic valuation of a sequence arising from
a rational integral, preprint.

\end{thebibliography}
\end{document}